# On the factorization of numbers of the form $X^2 + c$

July 22, 2022

**WOLF Marc,** https://orcid.org/0000-0002-6518-9882

Email: marc.wolf3@wanadoo.fr

**WOLF François,** https://orcid.org/0000-0002-3330-6087

Email: francois.wolf@dbmail.com

*Independent researchers;*

## Abstract

We study the factorization of the numbers $N = X^2 + c$, where $c$ is a fixed constant, and this independently of the value of $\gcd(X, c)$. We prove the existence of a family of sequences with arithmetic difference $(U_n, Z_n)$ generating factorizations, i.e. such that: $(U_n)^2 + c = Z_n Z_{n+1}$. The different properties demonstrated allow us to establish new factorization methods by a subset of prime numbers and to define a prime sieve. An algorithm is presented on this basis and leads to empirical results which suggest a positive answer to Landau's 4th problem.

# CONTENTS

# 1. Introduction

## 1.1. Primality and enumeration of prime numbers

Any odd integer $N$ can be written $N = X^2 + c$ where $X$ and $c$ are of opposite parity, $c \in \mathbb{N}^*$. When $c$ is fixed, we define $E_c = \{X^2 + c, X \in (2\mathbb{Z} + (c-1))\}$ the set of odd integers $N(X,c) \coloneqq X^2 + c$ thus defined. We are interested in the primes belonging to this set. Several methods exist to determine the primality of an odd number ([1], [2]) or more generally to determine all the primes of $[\![0, N]\!]$ ([3], [4], [5]). These methods cannot be easily extended to factorize an intermediate subset of numbers like $E_c \cap [\![0, N]\!]$ and thus to efficiently determine its primes when the size of the set becomes large.

In this paper, we first establish properties of the elements of $E_c$ and in particular of their factorization. We then establish that the elements of $E_c$ can be split into the product of two consecutive elements of a sequence with arithmetic difference $Z_n$ (not necessarily unique). This sequence is connected to a second sequence with arithmetic difference $U_n$ whose elements correspond to the values of $X$. Modular arithmetic then leads to a factorization of the integers of the set $E_c$ with a subset of prime numbers.

In the first part, we study the arithmetic sequences generated by odd prime factors which divide integers in the set $E_c$ and their powers. In the second part, we present pairs of sequences with arithmetic difference $(U_n, Z_n)$ which relate to the elements of $E_c$. Finally, we present an algorithm that factors the elements of $E_c$ and thus allows us to obtain all prime numbers of the form $X^2 + c$. This algorithm uses modular arithmetic in connection with the arithmetic sequences and the sequences with arithmetic difference previously introduced.

## 1.2. Definitions and notations

1- Let $\mathbb{P}$ be the set of prime numbers. Let $\mathcal{D}(x)$ be the (positive) divisors of a number x, and by extension let $\mathcal{D}(E)$ be the reunion of the divisors of all elements in $E$.
2- For any integer $(n, A) \in \mathbb{Z} \times \mathbb{N}^*$ let $\{n\}_A$ be the remainder of the Euclidean division of $n$ by $A$.
3- Let us define, for any $X \in \mathbb{Z}$, its <u>index</u> $j(X)$ as the quotient of its Euclidean division of $X$ by 2, i.e.:

$$X = 2j(X) + r, r \in \{0,1\}$$

The values $(X, c)$ considered will always be of opposite parity. For $c$ fixed, $r$ does not depend on $X$ and we will use this notation implicitly, i.e.:

$$r = 1 - \{c\}_2.$$

<u>Remark</u>: The indices of numbers were introduced in [6]. An application has been studied in [5].

4- For any $c \in \mathbb{N}^*$, there exists a unique pair $(t, y)$ such that:

$$c + 1 - r = 2^t(2y + 1)$$

5- For $A \in \mathcal{D}(E_c)$, we define $j_{A,c}(0)$ as the index of the smallest $X \coloneqq X_{A,c}(0)$ such that $A$ divides $N(X,c)$.

6- We define a sequence with arithmetic difference as a sequence $(U_n)$ such that $(U_{n+1} - U_n)$ is arithmetic.

# 2. Factorization and arithmetic sequences

The numbers in $E_c$ that have a common factor $A$ are connected to a set of arithmetic sequences. In this section, we present these sequences.

## 2.1. Arithmetic sequences and odd factors

**Proposition 2.1.1.** Let $A \in \mathcal{D}(E_c)$. For any element $N(X,c)$ of $E_c \cap A\mathbb{Z}$ there exists a divisor $a$ of $A$ and two natural numbers $j_{A,c,a}(0)$ and $k$ such that:

$$\gcd\left(a, \frac{A}{a}\right) | 2j_{A,c}(0) + r = X_{A,c}(0)$$

$$j(X) = j_{A,c,a}(0) + k.\frac{A}{\gcd\left(a, \frac{A}{a}\right)}$$

*Proof*: Let's assume that A divides $N(X,c)$, with $X = 2j(X) + r$, then we have:

$$N(X,c) = c + (2j(X) + r)^2 = AB$$

We know by definition of $j_{A,c}(0)$ that there exists $x \geq 0$ such that: $j(X) = j_{A,c}(0) + x$ and then we have:

$$\begin{aligned} N(X,c) &= c + \left(2j_{A,c}(0) + 2x + r\right)^2 \\ &= c + \left(2j_{A,c}(0) + r\right)^2 + 4x\left(2j_{A,c}(0) + r + x\right) \\ &= N\left(X_{A,c}(0), c\right) + 4x\left(2j_{A,c}(0) + r + x\right) \end{aligned}$$

Thus $A | x\left(2j_{A,c}(0) + r + x\right)$ hence there exists a factorization $A = ab$ such that $a|x$ and $b|2j_{A,c}(0) + r + x$, i.e. $x \equiv 0[a]$ and $x \equiv -2j_{A,c}(0) - r[b]$.

Invoking the Chinese theorem, we can find solutions if and only if $\gcd(a,b) | 2j_{A,c}(0) + r$.

Moreover, thanks to Bézout theorem, in this case we have a pair of integers $(u,v)$ such that $au - bv = -2j_{A,c}(0) - r$, so there must exist $k \in \mathbb{Z}$ such that $x = au + k.\frac{A}{gcd(a,b)}$. The choice of u minimal in $\mathbb{N}$ then yields $j_{A,c,a}(0) \coloneqq j_{A,c}(0) + au$ and $k \in \mathbb{N}$ as required.

Remark: The solutions of a similar linear Diophantine equation are given in [7].

**Definition 2.1.1.** Subject to the condition $\gcd\left(a, \frac{A}{a}\right) | 2j_{A,c}(0) + r$, we define the arithmetic sequence $\left(j_{A,c,a}(k)\right)$ by $j_{A,c,a}(k) \coloneqq j_{A,c,a}(0) + k.\frac{A}{\gcd\left(a, \frac{A}{a}\right)}$.

The sequences $\left(j_{A,c,a}(k)\right)$ and $\left(j_{A,c,\frac{A}{a}}(k)\right)$ exist under the same conditions, in which case they have the same common difference: we will call them dual sequences.

**Proposition 2.1.2.** Let $A \in \mathcal{D}(E_c)$. We have $j_{A,c}(0) \leq (A-r)/2$.

*Proof*: It is sufficient to show that there exists a natural number $X \leq A$ of parity opposite to that of c and such that $A|N(X,c)$.

We have $A|N(X,c) \Leftrightarrow X^2 \equiv -c[A]$. For any solution $x$ modulo $A$, $-x$ is also a solution. But if $A \in \mathcal{D}(E_c)$, then we can find such an $x$ between 0 and $A-1$, so that $x$ and $A-x$ are of opposite parity, we can therefore write $X = x$ or $X = A - x$ and conclude.

**Corollary 2.1.1.** With the notations of proposition 2.1.1, for $a \in \{1, A\}$, we have two dual sequences $\left(j_{A,c,1}(k)\right)$ and $\left(j_{A,c,A}(k)\right)$ with the same common difference $A$. These are the only two possible sequences when $A \in \mathbb{P}$. These two sequences are equal if and only if $A|c$. If they are distinct, we have the following relationship: $j_{A,c,1}(0) + j_{A,c,A}(0) = A - r$.

*Proof*: The condition $\gcd\left(a, \frac{A}{a}\right) | 2j_{A,c}(0) + r$ is trivially verified here since $\gcd\left(a, \frac{A}{a}\right) = 1$.

If $a = A$, we have $u = 0$ thus $j_{A,c,A}(0) = j_{A,c}(0)$.

If $a = 1$ we have $au \equiv -2j_{A,c}(0) - r[A]$.

The two sequences are equal if and only if $u = 0$ in both cases, i.e. $A|2j_{A,c}(0) + r = X_{A,c}(0)$. But we assumed $A|X_{A,c}(0)^2 + c$, which is equivalent to $A|c$ when $A$ is prime.

When on the contrary the two sequences differ, we must have $X_{A,c}(0) \neq 0$. Therefore, by minimality of $u$ and thanks to proposition 2.1.2, $au = A - 2j_{A,c}(0) - r$ and $j_{A,c,1}(0) = A - j_{A,c}(0) - r$. Thus, we have $j_{A,c,1}(0) + j_{A,c,A}(0) = A - r$.

**Remark 2.1.1.** More generally, the two dual sequences $\left(j_{A,c,a}(k)\right)$ and $\left(j_{A,c,\frac{A}{a}}(k)\right)$ are equal if and only if $\frac{A}{\gcd\left(a,\frac{A}{a}\right)} | X_{A,c}(0)$. $\frac{A}{\gcd\left(a,\frac{A}{a}\right)}$ is also the common difference for these two dual sequences.

## 2.2. Arithmetic sequences and prime power divisors

The numbers $N(X,c)$ in $E_c$ that share a prime factor $p$ are terms of one in at most two arithmetic sequences with the common difference equal to $p$. If we now look for the multiples of a power α strictly greater than 1 of $p$, i.e., $A = p^\alpha$, more sequences will appear. In this section, however, we show that one can always describe all terms of $E_c \cap A\mathbb{Z}$ with no more than two sequences, with common difference depending on $p, \alpha$ and $c$.

**Proposition 2.2.1:** Let $A = p^\alpha \in \mathcal{D}(E_c)$, $\gamma$ maximal such that $p^\gamma | X_{A,c}(0)$ (when $X_{A,c}(0) = 0$, $\gamma = +\infty$), and $\delta = \min\left(\gamma, \left\lfloor\frac{\alpha}{2}\right\rfloor\right)$. The sequence $\left(j_{A,c,a}(k)\right)$ exists if and only if $a \in \{p^\beta, p^{\alpha-\beta}\}_{\beta\leq\gamma}$. The elements of $E_c \cap A\mathbb{Z}$ are terms of one in at most two dual sequences $\left(j_{A,c,p^\delta}(k)\right)$ and $\left(j_{A,c,p^{\alpha-\delta}}(k)\right)$. These two sequences are equal if and only if $A|c$.

*Proof*: We use the same notations as in proposition 2.1.1. By hypothesis, any divisor $a$ of $A$ can be written $a = p^{\beta}, 0 \leq \beta \leq \alpha$.

The sequence $\left(j_{A,c,a}(k)\right)$ exists on the condition that gcd $\left(a, \frac{A}{a}\right) = p^{\min(\beta,\alpha-\beta)} | X_{A,c}(0)$ i.e. $\min(\beta, \alpha - \beta) \leq \gamma$. The common difference of the sequence is then equal to $p^{\max(\beta,\alpha-\beta)}$. In addition, we always have $\min(\beta, \alpha - \beta) \leq \frac{\alpha}{2}$, and conversely for any $\beta$ such that $0 \leq \min(\beta, \alpha - \beta) \leq \min\left(\gamma, \left\lfloor\frac{\alpha}{2}\right\rfloor\right) = \delta$ the sequence $\left(j_{A,c,p^{\beta}}(k)\right)$ exists and has common difference equal to $p^{\alpha-\beta}$.

Thus, $\delta$ is the maximum value that $\beta$ can take in the interval $\left[\!\left[0, \frac{\alpha}{2}\right]\!\right]$. We note that $\alpha - \delta \geq \delta$ thus so $\max(\delta, \alpha - \delta) = \alpha - \delta$.

For $\beta$ fixed, we distinguish the following two cases:

- If $\beta \geq \alpha - \beta$ (i.e. $\beta \geq \frac{\alpha}{2}$), we identify $u = x = 0$, because we have $x \equiv 0[p^{\beta}]$. The sequence $\left(j_{A,c,p^{\beta}}(k)\right)$ is then a subsequence of $\left(j_{A,c,p^{\alpha-\delta}}(k)\right)$, with common difference $p^{\beta} \geq p^{\alpha-\delta}$, and the same first term $j_{A,c}(0) = j_{A,c,p^{\alpha-\delta}}(0)$.
- Otherwise, when $\beta < \frac{\alpha}{2}$, $x$ is the smallest integer such that $x \equiv -X_{A,c}(0)[p^{\alpha-\beta}]$. The sequence $\left(j_{A,c,p^{\beta}}(k)\right)$ is again a subsequence of $\left(j_{A,c,p^{\delta}}(k)\right)$, with common difference $p^{\alpha-\beta} \geq p^{\alpha-\delta}$. Indeed, the first terms need not coincide but the congruence $x \equiv -X_{A,c}(0)[p^{\alpha-\beta}]$ implies $x \equiv -X_{A,c}(0)[p^{\alpha-\delta}]$ thus $j_{A,c,p^{\beta}}(0)$ is a term of the sequence $\left(j_{A,c,p^{\delta}}(k)\right)$.

If they exist, the sequences $\left(j_{A,c,p^{\beta}}(k)\right)$ and $\left(j_{A,c,p^{\alpha-\beta}}(k)\right)$ have the same common difference $p^{\max(\beta,\alpha-\beta)}$ and are therefore dual.

Finally, according to remark 2.1.1, the dual sequences $\left(j_{A,c,p^{\beta}}(k)\right)$ and $\left(j_{A,c,p^{\alpha-\beta}}(k)\right)$ are equal if and only if $p^{\max(\beta,\alpha-\beta)} | X_{A,c}(0)$, i.e. $\max(\beta, \alpha - \beta) \leq \gamma$. We deduce that we can distinguish two cases:

**Case 1)** $\boldsymbol{\gamma < \frac{\alpha}{2}}$, i.e. $A \nmid X_{A,c}(0)^2$: we get $\delta = \min\left(\gamma, \left\lfloor\frac{\alpha}{2}\right\rfloor\right) = \gamma$ and there exist $2(\gamma + 1)$ sequences, for $\beta \in \{0, \dots, \gamma, \alpha - \gamma, \dots, \alpha\}$. These sequences are pairwise distinct as per remark 2.1.1, because in this case their common difference never divides $X_{A,c}(0)$.

**Case 2)** $\boldsymbol{\frac{\alpha}{2} \leq \gamma}$, i.e. $A | X_{A,c}(0)^2$: The $\alpha + 1$ sequences are defined for $\beta \in \{0, \dots, \alpha\}$ and those with a common difference not greater than $p^{\gamma}$ are equal to their dual sequence.

As seen above, the dual sequences $\left(j_{A,c,p^{\delta}}(k)\right)$ and $\left(j_{A,c,p^{\alpha-\delta}}(k)\right)$ contain all the terms of the other sequences, i.e. any element of $E_c \cap A\mathbb{Z}$ is a term of one of these two sequences. They are equal only in case 2), when $\gamma \geq \left\lceil\frac{\alpha}{2}\right\rceil$, i.e. $p^{\alpha} | X^2$, or equivalently $A | c$.

**Corollary 2.2.1:** Let $A = p^{\alpha} \in \mathcal{D}(E_c), \alpha \in 2\mathbb{N}^*$. We let $X_{A,c}(0)^2 + c = AB$.

Then $p^\alpha | N(X,c)$ if and only if $j(X)$ is of one of the following:

- $j_{A,c,p^\delta}(k) = j_{A,c}(0) + p^\delta u + kp^{\alpha-\delta}$.

In this case, we have:

$$\begin{aligned} N(X,c) &= \left(2\left(j_{A,c}(0) + p^\delta u + kp^{\alpha-\delta}\right) + r\right)^2 + c \\ &= AB + 4\left(p^\delta u + kp^{\alpha-\delta}\right)\left(p^\delta u + kp^{\alpha-\delta} + X_{A,c}(0)\right) \\ &= p^\alpha \left[B + 4\left(u + kp^{\alpha-2\delta}\right)\left(\frac{X_{A,c}(0) + p^\delta u}{p^{\alpha-\delta}} + k\right)\right]. \end{aligned}$$

- $j_{A,c,p^{\alpha-\delta}}(k) = j_{A,c}(0) + kp^{\alpha-\delta}$.

In this case, we have:

$$\begin{aligned} N(X,c) &= \left(2\left(j_{A,c}(0) + kp^{\alpha-\delta}\right) + r\right)^2 + c \\ &= AB + 4kp^{\alpha-\delta}\left(kp^{\alpha-\delta} + X_{A,c}(0)\right) \\ &= p^\alpha \left[B + 4k\left(\frac{X_{A,c}(0)}{p^\delta} + kp^{\alpha-2\delta}\right)\right]. \end{aligned}$$

In the case where additionally $A \nmid c$, $p^{\alpha+1} \in \mathcal{D}(E_c)$ if and only if one of the following two equations (on $k$) has a solution:

$$B + 4u\left(\frac{X_{A,c}(0)+p^\gamma u}{p^{\alpha-\gamma}} + k\right) \equiv 0[p] \text{ or } B + 4k\left(\frac{X_{A,c}(0)}{p^\gamma}\right) \equiv 0[p].$$

*Proof*: The dichotomy on $j(X)$ is an immediate consequence of proposition 2.2.1. Furthermore if $A$ does not divide $c$, it implies that $\delta = \gamma < \frac{\alpha}{2}$, and the two Diophantine equations follow.

**Remark 2.2.1.** The previous corollary shows that, under the condition $p^\alpha \nmid c$, if $p^\alpha \in \mathcal{D}(E_c)$ then so are all powers of $p$.

*Proof*: Let us assume $p^\alpha \nmid c$ and $p^\alpha \in \mathcal{D}(E_c)$. Let us show by induction on $\beta \geq 0$ that $p^\beta \in \mathcal{D}(E_c)$. This is clearly true for $\beta$ from 0 to $\alpha$. If it is true for $\beta \geq \alpha$ fixed, the hypothesis $p^\alpha \nmid c$ implies $p^\beta \nmid c$ and therefore according to the above corollary, $p^{\beta+1} \in \mathcal{D}(E_c)$ if and only if one of the two equations has a solution:

$$B + 4u\left(\frac{X_{A,c}(0)+p^\gamma u}{p^{\alpha-\gamma}} + k\right) \equiv 0[p] \text{ or } B + 4k\left(\frac{X_{A,c}(0)}{p^\gamma}\right) \equiv 0[p].$$

These are Diophantine linear equations over $k$, so they have a solution as soon as their slope is invertible modulo $p$, i.e. non-divisible by $p$. The slopes are respectively $4u$ and $4\frac{X_{A,c}(0)}{p^\gamma}$. However, by maximality of $\gamma$, we know that $p \nmid \frac{X_{A,c}(0)}{p^\gamma}$.

Moreover $p$ being prime and odd, we deduce that $p \nmid 4\frac{X_{A,c}(0)}{p^\gamma}$.

**Remark 2.2.2.** Under the assumptions of the previous remark, the two linear Diophantine equations have a non-zero slope modulo p, in particular $u \not\equiv 0\ [p]$.

*Proof*: According to corollary 2.1.1, we have the following relation:

$$j_{p^\alpha,c,p^\gamma}(0) + j_{p^\alpha,c,p^{\alpha-\gamma}}(0) = p^\alpha - r.$$

Then we have $2j_{p^\alpha,c}(0) + p^\gamma u = p^\alpha - r$, hence:

$$\frac{X_{p^\alpha,c}(0)}{p^\gamma} = p^{\alpha-\gamma} - u.$$

Thus, we can deduce that $p \nmid u$ and conclude.

When we vary $j(X)$ in $\mathbb{N}$ and that $N(X,c) = X^2 + c$ is the first multiple of $p^\alpha$ found, the last proposition allows us to predict all subsequent occurrences of such multiples. This can be used in the algorithm presented in the last part.

**Proposition 2.2.2.** We keep the assumptions and notations of proposition 2.2.1. If $p \nmid c$, the common difference of the two dual sequences $\left(j_{A,c,p^\delta}(k)\right)$ and $\left(j_{A,c,p^{\alpha-\delta}}(k)\right)$ is $p^\alpha$, and the sum of their initial values is $p^\alpha - r$. Otherwise, let $\nu$ be maximal such that $p^\nu | c, \nu \in \mathbb{N}^*$. If $\alpha \le \nu$, a single self-dual sequence can describe all the elements of $E_c \cap p^\alpha\mathbb{Z}$, with common difference $p^{\lceil \alpha/2 \rceil}$. If $\alpha > \nu$, $p^\alpha \in \mathcal{D}(E_c)$ if and only if $\nu$ is even and $p \in \mathcal{D}\left(E_{c/p^\nu}\right)$. We can then describe the elements of $E_c \cap p^\alpha\mathbb{Z}$ with two dual sequences with common difference $p^{\frac{\nu}{2}+(\alpha-\nu)}$, which initial values sum to $p^{\frac{\nu}{2}+(\alpha-\nu)} - r$.

*Proof*: If $p \nmid c$ then $p \nmid X_{A,c}(0)$ thus $\delta = 0$. In particular according to proposition 2.2.1 the two dual sequences $\left(j_{A,c,1}(k)\right)$ and $\left(j_{A,c,A}(k)\right)$ are distinct and the corollary 2.1.1 ensures that $j_{A,c,1}(0) + j_{A,c,A}(0) = A - r$.

Now let $\nu$ be maximal such that $p^\nu | c$.

If $\alpha \le \nu$, letting $X = 0$ or $X = p^\alpha$ we observe that $p^\alpha \in \mathcal{D}(E_c)$ and the existence of the announced self-dual sequence is a direct consequence of proposition 2.2.1.

If $\alpha > \nu$ and $\nu$ is odd, then $p^\alpha \in \mathcal{D}(E_c)$ implies $p^\nu | X_{p^\alpha,c}(0)^2$, and the exponent of $p$ in the decomposition of $X_{p^\alpha,c}(0)^2$ into prime factors being necessarily even, we necessarily have:

$$p^{\nu+1} | X_{p^\alpha,c}(0)^2$$

thus $p^{\nu+1} | c$ which is impossible.

Finally if $\alpha > \nu$ and $\nu$ even, we have $p^\alpha | X^2 + c$ if and only if $p^{\frac{\nu}{2}} | X$ and $p^{\alpha-\nu} \left| \left(\frac{X}{p^{\frac{\nu}{2}}}\right)^2 + \frac{c}{p^\nu}\right.$.

Thus, the elements of $\mathcal{D}(E_c) \cap p^\alpha\mathbb{Z}$ correspond one-to-one to those of $\mathcal{D}\left(E_{c/p^\nu}\right)$. However, we know that $c/p^\nu$ is prime with $p$, thus the final expected result follows from remark 2.2.1 and the first result of this proposition.

## 2.3. Asymptotics of the number of prime factors

Let $E_c^{(J)} = E_c \cap [\![0,(2J+r)^2+c]\!]$. Clearly $\left|E_c^{(J)}\right| = J+1$.

The sets $\mathcal{D}\left(E_c^{(J)}\right)$ and $\mathcal{D}\left(E_c^{(J)}\right) \cap \mathbb{P}$ increase with $J$.

We show in this section that there is a threshold value $J_c$ from which $\mathcal{D}\left(E_c^{(J)}\right) \cap \mathbb{P}$ grows by at most one prime number with each increment of $J$.

**Proposition 2.3.1.** Let $p \in \mathbb{P}$, $\alpha \in \mathbb{N}^*$. If $p^\alpha \in \mathcal{D}(E_c)$, then $j_{p^\alpha,c}(0) \le (p^\alpha - 1)/2$, with equality if and only if one of the following Is true:

- $c$ is even and $p^\alpha | c$
- $c$ is odd and $p^\alpha | c+1$.

*Proof*: This is a refinement of proposition 2.1.2. Let us assume $p^\alpha \in \mathcal{D}(E_c)$. If $p^\alpha | c$ then $j_{p^\alpha,c}(0) = 0$ or $\frac{p^\alpha-1}{2}$ depending on the value of $r$. Otherwise, according to corollary 2.1.1, we have two distinct dual sequences which initial values sum to $p^\alpha - r$, thus $j_{p^\alpha,c}(0) < \frac{p^\alpha - r}{2}$ and $j_{p^\alpha,c}(0) \le \frac{p^\alpha-1}{2}$.

**Proposition 2.3.2.** Let $c > 0$. Let:

$$J_c = \left\lfloor \frac{c-1}{4} \right\rfloor$$

For any $J > J_c$ we have:

$$\left|\mathcal{D}\left(E_c^{(J)}\right) \cap \mathbb{P}\right| - \left|\mathcal{D}\left(E_c^{(J-1)}\right) \cap \mathbb{P}\right| \le 1.$$

In addition, if $p \in \left(\mathcal{D}\left(E_c^{(J)}\right) \cap \mathbb{P}\right) \setminus \mathcal{D}\left(E_c^{(J-1)}\right)$ then $p^2 \notin \mathcal{D}\left(E_c^{(J)}\right)$.

*Proof*: We assume that $J$ is such that $\left|\mathcal{D}\left(E_c^{(J)}\right) \cap \mathbb{P}\right| - \left|\mathcal{D}\left(E_c^{(J-1)}\right) \cap \mathbb{P}\right| \ge 2$, or $p^2 \in \left(\mathcal{D}\left(E_c^{(J)}\right) \cap \mathbb{P}\right) \setminus \left(\mathcal{D}\left(E_c^{(J-1)}\right) \cap \mathbb{P}\right)$.

There are therefore $p < q \in \mathbb{P}$ such that $j_{p,c}(0) = j_{q,c}(0) = J$, or $j_{p,c}(0) = j_{p^2,c}(0) = J$.

In particular according to the proposition 2.1.2 we have $J \le (p-1)/2$. Equality is only possible if $p|c$ with $c$ even or $p|c+1$ with $c$ odd.

1. Let us first assume that we are not in one of these two cases of equality, i.e. $J \le (p-3)/2$. Then with $X = 2J + r$ we have $X^2 + c$ that is multiple of $pq$ or $p^2$ and therefore, in both cases, it is greater than $p^2$ which gives:

$$p^2 \le (2J+r)^2 + c \le (p+r-3)^2 + c$$

Thus $c \ge p^2 - p^2 + 2(3-r)p - (3-r)^2$.

So $c \ge 2(3-r)p - (3-r)^2 \ge 2(3-r)(2J+3) - (3-r)^2 = (3-r)(4J+3+r)$.

In other words, we get $J \le \left\lfloor \frac{c-9+r}{4(3-r)} \right\rfloor$.

2. <u>If $c$ is odd and $p$ divides $c+1$</u>, we have $J \leq (p-1)/2$ and $X^2 + c \geq (p-1)^2 + c$. We then distinguish two cases:
   a) If $j_{q,c}(0) = J$, we have $q \geq p+2$ and so:

$$p(p+2) \leq (p-1)^2 + c = p^2 - 2p + 1 + c$$

hence $4p \leq c+1$ and $J = \frac{p-1}{2} \leq \left\lfloor \frac{c-3}{8} \right\rfloor$.

   b) If $j_{p^2,c}(0) = J$, we have $X^2 + c = (p-1)^2 + c \geq p^2$ hence $c + 1 \geq 2p$. Thus $J = \frac{p-1}{2} \leq \left\lfloor \frac{c-1}{4} \right\rfloor$.

<u>Example</u>: for $c = 61$, we then have $J = \left\lfloor \frac{c-1}{4} \right\rfloor = 15$ with $N = 30^2 + 61 = 31^2$.

3. <u>If $c$ is even and $p$ divides $c$</u>, we have $J \leq (p-1)/2$ et $X^2 + c \geq p^2 + c$. We then distinguish two cases:
   a) If $j_{q,c}(0) = J$, we have $q \geq p+2$ and thus:

$$p(p+2) \leq p^2 + c$$

hence $c \geq 2p$ and $J = \frac{p-1}{2} \leq \left\lfloor \frac{c-2}{4} \right\rfloor$.

   b) If $j_{p^2,c}(0) = J$, we have $X^2 + c > p^2$ and $p^2 | X^2 + c$ therefore $X^2 + c \geq 3p^2$. Thus $c \geq 2p^2$ and $J \leq \left\lfloor \sqrt{\frac{c}{8}} - \frac{1}{2} \right\rfloor$.

We deduce that, when c is even, $J \leq \max\left(\left\lfloor \frac{c-8}{8} \right\rfloor, \left\lfloor \sqrt{\frac{c}{8}} - \frac{1}{2} \right\rfloor, \left\lfloor \frac{c-2}{4} \right\rfloor, 0\right) = \left\lfloor \frac{c-2}{4} \right\rfloor$, and when $c$ is odd, $J \leq \max\left(\left\lfloor \frac{c-9}{12} \right\rfloor, \left\lfloor \frac{c-3}{8} \right\rfloor, \left\lfloor \frac{c-1}{4} \right\rfloor, 0\right) = \left\lfloor \frac{c-1}{4} \right\rfloor$.

In all cases, $J > J_c$ implies the expected results on the appearance of new prime factors in $N(2J + r, c)$.

**<u>Remark 2.3.1.</u>** For any $X > 2J_c + r$, we have: $X = \left\lfloor \sqrt{c + X^2} \right\rfloor$.

*<u>Proof</u>*: Clearly $\sqrt{c + X^2} \geq X$. Therefore $X = \left\lfloor \sqrt{c + X^2} \right\rfloor$ if and only if $\sqrt{c + X^2} < X + 1$ thus $c < 2X + 1$ and with $X = 2j + r$, we get:

$$j > \frac{c - (2r+1)}{4}$$

Thus if $j > J_c = \left\lfloor \frac{c-1}{4} \right\rfloor$, the inequality is verified.

**<u>Proposition 2.3.3.</u>** $\mathcal{D}(E_c) \cap \mathbb{P}$ is infinite.

*<u>Proof</u>*: We assume that $\mathcal{D}(E_c) \cap \mathbb{P}$ is finite, necessarily non-empty. We let $p_1 \ldots p_m$ be its elements prime with $c$ and $q_1 \ldots q_n$ those that divide $c$. We set:

$$X = 2^{1-r}\left(\prod p_i\right)\left(1 + \prod q_j\right).$$

We note that $X > 1$. We have then, for any prime number $p_i$:

$$X^2 + c \equiv c[p_i]$$

And for any prime number $q_j$:

$$X^2 + c \equiv \left(2^{1-r}\left(\prod p_i\right)\right)^2 [q_j]$$

Since $c \notin p_i\mathbb{Z}$ and $2^{1-r}(\prod p_i) \notin q_j\mathbb{Z}$ we deduce that $X^2 + c$ has a prime factor distinct from $p_i$ and $q_j$, which is impossible.

**Proposition 2.3.4.** $\mathcal{D}(E_c) \cap \mathbb{P}$ is smaller than $\mathbb{P}$.

*Proof*: We will show the stronger following asymptotic result:

$$\frac{\left|\mathcal{D}\left(E_c^{(J)}\right) \cap \mathbb{P}\right|}{|[\![0, N(2J + r, c)]\!] \cap \mathbb{P}|} \to_{J\to+\infty} 0.$$

When $J \to +\infty$, the prime number theorem ensures that:

$$|[\![0, N(2J + r, c)]\!] \cap \mathbb{P}| \sim \frac{N(2J + r, c)}{\ln\big(N(2J + r, c)\big)} \sim \frac{2J^2}{\ln J}$$

On the other hand, we have shown in proposition 2.3.2 that beyond the threshold $J_c$, we have:

$$\left|\mathcal{D}\left(E_c^{(J)}\right) \cap \mathbb{P}\right| - \left|\mathcal{D}\left(E_c^{(J-1)}\right) \cap \mathbb{P}\right| \leq 1$$

We can therefore deduce that $\left|\mathcal{D}\left(E_c^{(J)}\right) \cap \mathbb{P}\right| = O(J)$.

Thus:

$$\frac{\left|\mathcal{D}\left(E_c^{(J)}\right) \cap \mathbb{P}\right|}{|[\![0, N(2J + r, c)]\!] \cap \mathbb{P}|} = O\left(\frac{\ln(J)}{J}\right) \to_{J\to+\infty} 0.$$

# 3. Pairs of sequences that generate factorizations

## 3.1. Sequences with arithmetic difference

In this section, we will exhibit pairs of sequences with arithmetic difference that generate factorizations of elements $N(X, c)$ of $E_c$. These sequences $(U_n, Z_n)$ verify the following relation:

$$(1) \qquad U_n^2 + c = Z_n Z_{n+1}$$

Consider the general case of a sequence with arithmetic difference $(U_n)_{n\in\mathbb{Z}}$ defined by:

$$U_n = \alpha n(n - 1) + \beta n + U_0.$$

To be able to use $U_n$ in $N(X, c)$ we impose that $\alpha \in \mathbb{Z}$, $\beta \in 2\mathbb{Z}$ and $U_0 \in 2\mathbb{Z} + (c - 1)$, so that all the terms of the sequence $(U_n)$ have the same parity, opposite to $c$.

For such a sequence, we set $Z_n = \frac{1}{2}(U_n + U_{n-1}) = U_n - \alpha(n - 1) - \frac{\beta}{2} \in \mathbb{Z}$. Thus, we have:

$Z_n = \alpha n(n-1) + \gamma n + Z_0$, with $\gamma = \beta - \alpha$ and $Z_0 = \alpha - \frac{\beta}{2} + U_0$.

**Proposition 3.1.1.** The sequences defined above verify (1) when $c = \alpha U_0 + \frac{\alpha\beta}{2} - \frac{\beta^2}{4}$.

*Proof:* An elementary calculation shows that:

$$\begin{aligned} Z_n Z_{n+1} &= \left(U_n - \alpha(n-1) - \frac{\beta}{2}\right)\left(U_n + \alpha n + \frac{\beta}{2}\right) \\ &= U_n^2 + \alpha U_n - \left(\alpha^2 n(n-1) + \alpha\beta n - \frac{\alpha\beta}{2} + \frac{\beta^2}{4}\right) \\ &= U_n^2 + \left(\alpha U_0 + \frac{\alpha\beta}{2} - \frac{\beta^2}{4}\right) \end{aligned}$$

We can then conclude.

**Proposition 3.1.2.** For any $N(X,c) \in E_c$ and $A|N(X,c)$, there exist a unique pair $(U_n, Z_n)$ verifying (1) and such that $U_0 = X, Z_0 = A$.

*Proof:* We write $N(X,c) = X^2 + c = AB$. Keeping the same notations, we have $Z_0 = \alpha - \frac{\beta}{2} + U_0$. Therefore, we want to solve the system:

$$\begin{cases} \alpha - \dfrac{\beta}{2} + X = A \\ \alpha X + \dfrac{\alpha\beta}{2} - \dfrac{\beta^2}{4} = c \end{cases} \quad (\alpha, \beta) \in \mathbb{Z} \times 2\mathbb{Z}$$

Equivalently:

$$\begin{cases} \dfrac{\beta}{2} - \alpha = X - A \\ \alpha X - \dfrac{\beta}{2}\left(\dfrac{\beta}{2} - \alpha\right) = c \end{cases}$$

The second equation gives $X\alpha - \frac{\beta}{2}(X - A) = c$ then $X\left(\alpha - \frac{\beta}{2}\right) + A\frac{\beta}{2} = c$ thus:

$$\frac{\beta}{2} = \frac{1}{A}(X^2 + c - AX) = B - X.$$

We hence get:

$$\begin{cases} \beta = 2(B - X) \\ \alpha = A + B - 2X \end{cases}$$

Note that we also have: $\gamma = \beta - \alpha = B - A$.

*Remark*: We do not necessarily have $A|U_n^2 + c$ for the other terms of the sequence.

On the other hand, if we let $U_n(X, A, c)$ be the sequence thus obtained, we can describe all possible factorizations of the integers $N(X,c)$ by changing $X$ and $A$. It would be interesting to try and determine a subfamily of elements$(X, A)$ sufficient to generate all possible factorizations. The asymptotic density of the integers of the sequence $(U_n)$ being zero, such a subfamily would necessarily be infinite.

## 3.2. Arithmetic sequences and sequences with arithmetic difference

For any $N(X,c) = c + (2j(X) + r)^2$ of $E_c \cap A\mathbb{Z}$, the arithmetic sequence with common difference $A$ of indices $j_{A,c}^X(k) \coloneqq \{j(X)\}_A + k.A$ corresponds to elements of $E_c \cap A\mathbb{Z}$ including $N(X,c)$. This sequence corresponds (proposition 3.1.2) to a family of pairs of sequences with arithmetic difference $\big(U_n(k), Z_n(k)\big)_{n\in\mathbb{Z}}$ such that, for any $n$ and $k$:

$$Z_0 = A$$

$$U_0(k) = 2j_{A,c}^X(k) + r$$

$$U_n^2(k) + c = Z_n(k)Z_{n+1}(k)$$

In this section, we explicit the coefficients $\alpha(k)$, $\beta(k)$ et $\gamma(k)$ of these sequences.

**Proposition 3.2.1.** With the previous notations, if we additionally let $X_{A,c}^X(k) = 2j_{A,c}^X(0) + r + 2kA$ and $X_{A,c}^X(0)^2 + c = AB$, we have:

$$U_n(k) = \alpha(k)n(n-1) + \beta(k)n + X_{A,c}^X(k)$$

$$Z_n(k) = \alpha(k)n(n-1) + \gamma(k)n + A$$

with:

$$\begin{cases} \alpha(k) = 4Ak(k-1) + 4X_{A,c}^X(0)k + \left(A + B - 2X_{A,c}^X(0)\right) \\ \beta(k) = 8Ak(k-1) + \left(8X_{A,c}^X(0) + 4A\right)k + 2\left(B - X_{A,c}^X(0)\right) \\ \gamma(k) = 4Ak(k-1) + \left(4X_{A,c}^X(0) + 4A\right)k + (B - A) \end{cases}$$

These coefficients are sequences with arithmetic difference over k.

*Proof:* This is an immediate consequence of proposition 3.1.2, which can also be verified directly by calculation.

Remark: More generally, any arithmetic sequence of indices introduced in part 2 can be mapped to a sequence of pairs of sequences with arithmetic difference, the initial term of which corresponds to a trivial factorization of an element of $E_c$. The coefficients of these sequences are given in Appendix A.

## 3.3. Two particular sequences

A first special case is used in the factorization algorithm, stemming from the following factorization (see definition 4 in section 1.2):

$$(1-r)^2 + c = c + 1 - r = 2^t(2y+1)$$

which corresponds to a factorization of an even number $X^2 + c$ by a minimal odd factor. To avoid the parity issue, we let $X = (-1)^r(2y + r)$, so that $X \equiv 1 - r[2y+1]$ and the factorization $N(X,c) = (2y+r)^2 + 2^t(2y+1) - (1-r)$ by $A = 2y + 1$ remains possible. Indeed, we have $N(X,c) \equiv (1-r)^2 - (1-r) \equiv 0\ [2y+1]$ and we can even show that $B = (2y+1) - 2(1-r) + 2^t$.

We then simply determine the unique pair given by proposition 3.1.2:

$$\begin{cases} \beta = 2[(2y+1) - (-1)^r(2y+r) - 2(1-r) + 2^t] \\ \alpha = 2[(2y+r) - (-1)^r(2y+r) + 2^{t-1}] \end{cases}$$

By substituting $(-1)^r$ by $(1-2r)$, we obtain a simpler expression:

$$\begin{cases}\beta = 2\left[2^t - 1 + r\big(1 + 2(2y+1)\big)\right] \\ \alpha = 2[2r(2y+1) + 2^{t-1}]\end{cases}$$

We recall that $Z_n = \alpha n(n-1) + (\beta - \alpha)n + A$.

**Proposition 3.3.1.** We set $Z_n$ as above. The values $Z_n$ are pairwise distinct for $n \in \mathbb{Z}$ if and only if $t > 1$ or $r = 1$.

*Proof*: The points $(n, Z_n)$ are on the parabola $y = \alpha x^2 + (\beta - 2\alpha)x + A$, which is minimal for $x_{min} = -\left(\frac{\beta}{2\alpha} - 1\right)$. For $n \neq m$, $Z_n = Z_m$ if and only if $n + m = 2x_{min} = -\left(\frac{\beta}{\alpha} - 2\right)$ therefore the values $Z_n$ are pairwise distinct for $n \in \mathbb{Z}$ if and only if $\alpha \nmid \beta$. Furthermore:

$$\frac{\beta}{\alpha} = \frac{2^t - 1 + r\big(1 + 2(2y+1)\big)}{2r(2y+1) + 2^{t-1}}$$

If $r = 0$, $\alpha|\beta$ is equivalent to $2^{t-1}|2^t - 1$ i.e. $t = 1$.

If $r = 1$, $\alpha|\beta$ is equivalent to $2^{t-1} + 2(2y+1)|2^t + 2(2y+1) = 2.2^{t-1} + 2(2y+1)$ i.e. $2^{t-1} + 2(2y+1)|2^{t-1}$ which is impossible.

We have thus established the expected result.

We now introduce a second pair $(U_n, Z_n)$, defined only if $t \geq 4 - r$, from the following factorization:

$$(3-r)^2 + c = c + 9 - 5r = 2^t(2y+1) + 8 - 4r$$

By hypothesis, $2^t(2y+1) + 8 - 4r = 2^{3-r}(2^{t-3+r}(2y+1) + 1)$ and $2z + 1 = 2^{t-3+r}(2y+1) + 1$ is odd. As before, we avoid the parity issue ($(3-r)^2 + c$ is even) by setting $X = -(3-r) + (3-2r)(2z+1)$, so that $X^2 + c$ is a multiple of $2z+1$, more precisely, using that $9 - 5r = 2^{3-r} + 1 - r$, we have $X^2 + c = A(2z+1)$ with $A = \frac{1}{2^{3-r}}(c + 2^{3-r} + 1 - r)$.

*Note:* the choice of the $3 - 2r$ factor rather than 1 is arbitrary and could be changed.

The proof of the following expressions is left to the reader:

**Case 1)** $r = 0$: with $U_0 = 3(c_1+1)/8$ and $Z_0 = 1 + (c_1+1)/8$, and:

$$U_n = \frac{(c_1+1)}{2}n(n-1) + \frac{(3c_1-1)}{2}n + U_0$$

$$Z_n = \frac{(c_1+1)}{2}n(n-1) + (c_1-1)n + Z_0$$

**Case 2)** $r = 1$: with $U_0 = c_2/4 - 1$ and $Z_0 = c_2/4 + 1$, and:

$$U_n = 4n(n-1) + 4n + U_0$$

$$Z_n = 4n(n-1) + Z_0$$

**Proposition 3.3.2.** We set $(U_n, Z_n)$ as above. If $r = 0$, the values of $Z_n$ for $n \in \mathbb{Z}$ are pairwise distinct.

*Proof*: As in proposition 3.3.1, this amounts to proving that $\frac{\beta}{\alpha}$ is not an integer, however $\frac{\beta}{\alpha} = \frac{3c-1}{c+1} = 3 - \frac{4}{c+1}$. As we have assumed $16|c+1$ we deduce the expected result.

# 4. An Algorithm to factorize integers of the set $E_c$ and related prime sieve

## 4.1. Presentation of the algorithm

We present here an algorithm which aims at finding the prime numbers in the set $E_c^{(J)}$. Most of the results of the previous sections allow us to crucially speed up the factorization of the elements of $E_c$. The pairs of sequences with arithmetic difference from part 3 allow us in particular to generate "free" factorizations for small values of $j$.

We note $P_J = \mathbb{P} \cap E_c^{(J)}$ and $D_J = \mathcal{D}\left(E_c^{(J)}\right) \cap \mathbb{P} \cap [\![0, J]\!]$. The algorithm returns these two sets.

**ALGORITHM**:

For the sake of brevity, the detailed algorithm has been implemented for $r = 0$. This algorithm is available online at the following URL:
https://drive.google.com/file/d/1qnF9bfYSBCxna8eDmi5j85tExmqRQzy2/view?usp=sharing

We give below a general description of the algorithm independently of the value of $r$:

- The inputs of the algorithm are the numbers $c$ and $J$.
- The first step is to initialize variables such as $J_c$ and arrays for prime number and factorization storage.
- The second step is to determine the prime numbers in $[\![0, \sqrt{N(2J_c + r, c)}\,]\!]$ by using Atkin's algorithm. This algorithm uses little memory, especially when number indices are used. More details can be found in [5] on this subject. In this step, we also determine the first prime factors using the two propositions of section 3.3 and the previous array of primes. The integer $c + r$ is the first odd number factorized.
- The third step performs the factorization of all numbers in $E_c^{(J)}$, partitioned in two subsets: first the numbers not greater than $N(2J_c + r, c)$ for which we complete the factorizations initiated in the previous step, and then the rest, for which we can purely rely on arithmetic progressions to predict factorizations and easily pick new primes as they appear at most one at a time.
- As full factorizations are performed, we can also calculate sets $P_J$ and $D_J$. It is also possible, for verification purposes, to get the algorithm to additionally return the details of all factorizations performed.

The algorithm thus determines the primality of the elements in $E_c$ by performing a sieve-type method; that is, after a certain rank, such an element is determined to be prime when it is not a term of any of the sequences generated by the previously determined factors. However, the complete factorization of the composite elements is necessary to accurately detect the new prime divisors that are not in $E_c$ and keep the algorithm going.

Remark: when $c \leq 3$, the first subset is empty. For $J_c = 0$ we have only $N(r,c) = r + c$ to factor in the second step.

## 4.2. Results

In this section only, we have set $c = 1$.

We will first present the empirical complexity (computation time) of our algorithm. We will then compare the cardinality of the sets $P_J$ and $D_J$. Details of the Maple options used to obtain the linear regressions presented in this section are given in Appendix B.

On figure 1 below, we represent the execution time of the algorithm, in seconds, as a function of $N(2J,1)$ up to a value $J = 10^8$, that is $N(2J,1) = 4.10^{16} + 1$. We obtain an empirical complexity in $O(N(2J,1))$.

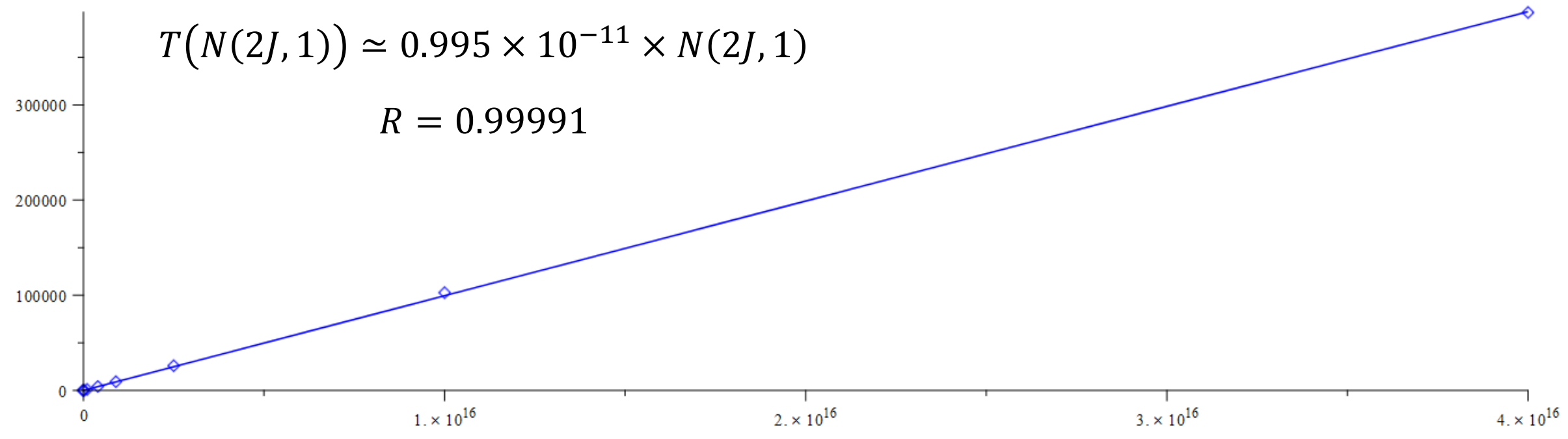

**Figure 1:** calculation time (s) for $c = 1$, function of $N(2J,1)$

In [5], we compared the Atkin sieve method and the Pritchard wheel sieve method. The complexity in terms of number of operations of these two methods is $O(N)$, i.e. empirically the same as the algorithm presented in this paper (the former to determine all primes, the latter only a subset). However, in practice the latter achieves prime numbers almost 400 times larger in the same time of the order of a second.

On figure 2 below, we represent $|P_J|$ and $|D_J|$ as a function of $J$. When $J \to +\infty$, proposition 2.3.3 ensures that $|D_J| \to +\infty$. Empirically, we can conjecture that the same is true for $|P_J|$.

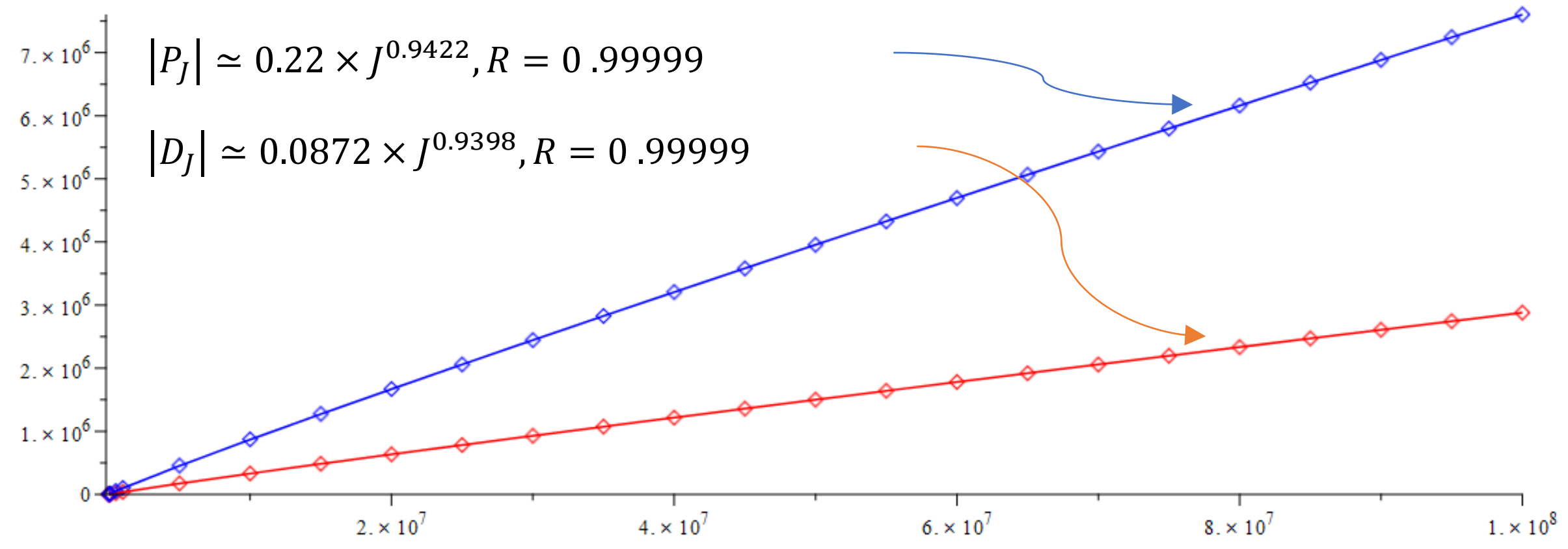

**Figure 2:** $|P_J|$ and $|D_J|$ function of $J$

This empirical result suggests a positive answer to the fourth Landau problem. Similar results having been obtained for the values of $c > 1$ studied, thus we propose the following:

Conjecture: For any $c \in \mathbb{N}^*$, there exists an infinity of prime numbers of the form $X^2 + c$.

## 5. Conclusion

We have studied the factorization of integers of the form $N(X, c) = X^2 + c$ via dual arithmetic sequences on $X$. We have given several asymptotic results on the occurrence of prime factors in $N(X, c)$. We have also introduced a family of pairs of sequences with arithmetic difference $(U_n, Z_n)$ which provide the following factorization:

$$U_n^2 + c = Z_n . Z_{n+1}$$

We have applied the conjunction of these results to implement a prime sieve on $E_c$, and presented our results for $c = 1$. We subsequently established a link with Landau's fourth problem, without solving it categorically; more generally we conjectured that for any $c \in \mathbb{N}^*$, there exists an infinity of primes of the form $X^2 + c$.

Acknowledgments: We would like to thank François-Xavier VILLEMIN for his attentive comments and suggestions.

## Appendix A: Presentation of sequences with two parameters k and j

Let $N = N(X,c) = X^2 + c$ be an element of $E_c$, and $X = 2j + r$. The factorization of $N = N \times 1$ generates, as given in part 2.1, two dual arithmetic sequences of indices which are multiples of $N$ in $E_c$:

$$\begin{cases} j(k) = N - j - r + kN \\ j'(k) = j + (k+1)N \end{cases}$$

Each factorization of the elements of $E_c$ corresponding to the above indices gives rise to a pair of sequences with arithmetic difference as shown in proposition 3.1.2. The first terms are naturally given by:

$$\begin{cases} U_0(k,j) = 2[(N - j - r) + kN] + r & (\boldsymbol{A}) \\ U_0'(k,j) = 2[j + (k+1)N] + r & (\boldsymbol{B}) \end{cases}$$

The pair of sequences with arithmetic difference $\left(U_n(k,j), Z_n(k,j)\right)$ obtained by using the first term $(\boldsymbol{A})$ is written as follows:

$$\begin{cases} U_n(k,j) = \left(\alpha_1(k) + \alpha_2(k,j)\right)n(n-1) + \left(\beta_1(k) + \beta_2(k,j)\right)n + U_0(k,j) \\ Z_n(k,j) = \left(\alpha_1(k) + \alpha_2(k,j)\right)n(n-1) + \left(\gamma_1(k) + \gamma_2(k,j)\right)n + Z_0(j) \end{cases}$$

where the coefficients $\alpha_1(k)$, $\beta_1(k)$ and $\gamma_1(k)$ are given by:

$$\begin{cases} \alpha_1(k) = 4ck(k+1) + \left(c + (1-r)\right) \\ \beta_1(k) = 4ck(k+1) + cV_k + 2\left(c + (1-r)\right) \\ \gamma_1(k) = cV_k + \left(c + (1-r)\right) \end{cases}$$

the sequence $V_k = 4k(k-1) + 12k + 2$, and:

$$\begin{cases} \alpha_2(k,j) = 4j\left(j - (1-r)\right) + 24jk + 32kj(j-1) + 16j(j-r)k(k-1) \\ \qquad\qquad +4kr(8j + 8j(k-1) + k) \\ \beta_2(k,j) = 4j + 16j\left(j - (1-r)\right) + 64jk + 80kj(j-1) + 32j(j-r)k(k-1) \\ \qquad\qquad +4kr(20j + 16j(k-1) + 2k + 1) \\ \gamma_2(k,j) = 4j + 12j\left(j - (1-r)\right) + 40jk + 48kj(j-1) + 16j(j-r)k(k-1) \\ \qquad\qquad +4kr(12j + 8j(k-1) + k + 1) \end{cases}$$

The pair of sequences with arithmetic difference $\left(U_n'(k,j), Z_n'(k,j)\right)$ obtained by using the first term $(\boldsymbol{B})$ is written in a similar way to the previous pair of sequences with the coefficients $\alpha_1'(k)$, $\beta_1'(k)$ and $\gamma_1'(k)$ given by:

$$\begin{cases} \alpha_1'(k) = \alpha_1(k) + \boldsymbol{4r} \\ \beta_1'(k) = \beta_1(k) + \boldsymbol{12r} \\ \gamma_1'(k) = \gamma_1(k) + \boldsymbol{8r} \end{cases}$$

And:

$$\begin{cases} \alpha_2'(k,j) = \alpha_2(k,j) + \boldsymbol{8(j + 2jk) + 8kr} \\ \beta_2'(k,j) = \beta_2(k,j) + \boldsymbol{8(3j + 4jk) + 16kr} \\ \gamma_2'(k,j) = \gamma_2(k,j) + \boldsymbol{8(2j + 2jk) + 8kr} \end{cases}$$

The families of pairs of sequences with arithmetic difference thus obtained allow us to generate a large number of factorizations (but not necessarily all).

## Appendix B: Maple regressions

Here are numeric values obtained from our implementation (Visual Studio C++ 2019) of the algorithm presented in this article.

Table 1 contains numeric values of $T\big(N(2J,1)\big)$ (in seconds) measured from running the algorithm with $c = 1$.

**Table 1: numeric values of $T\big(N(2J,1)\big)$ in seconds.**

| $J$ | $1 \times 10^5$ | $5 \times 10^5$ | $1 \times 10^6$ | $5 \times 10^6$ | $1 \times 10^7$ | $1{,}5 \times 10^7$ |
|---|---|---|---|---|---|---|
| $T\big(N(2J,1)\big)$ | 0,592 | 12,323 | 47,147 | 1047,542 | 4099,852 | 9188,870 |

| $J$ | $2{,}5 \times 10^7$ | $5 \times 10^7$ | $1 \times 10^8$ |
|---|---|---|---|
| $T\big(N(2J,1)\big)$ | 25957,098 | 102681,071 | 397131,895 |

For the linear regression of $T\big(N(2J,1)\big)$, we used Maple's routine **Fit** as below:

$$\texttt{Fit}(a \times n + b,\ \texttt{X},\ \texttt{Y},\ n,\ \texttt{summarize = embed})$$

We get the following regression with $b$ set to 0:

$$T\big(N(2J,1)\big) \simeq 9.95003248816399 \times 10^{-12} \times N(2J,1),\ R = 0.999908$$

Table 2 contains the values of $|P_J|$ and $|D_J|$ returned by the algorithm with $c = 1$.

**Table 2: values of $|P_J|$ and $|D_J|$.**

| $J$ | $1 \times 10^4$ | $5 \times 10^4$ | $1 \times 10^5$ | $5 \times 10^5$ | $1 \times 10^6$ | $5 \times 10^6$ |
|---|---|---|---|---|---|---|
| $\lvert D_J \rvert$ | 609 | 2549 | 4783 | 20731 | 39175 | 174193 |
| $\lvert P_J \rvert$ | 1558 | 6655 | 12390 | 54109 | 102204 | 456361 |

| $J$ | $1 \times 10^7$ | $1,5 \times 10^7$ | $2 \times 10^7$ | $2,5 \times 10^7$ | $3 \times 10^7$ |
|---|---|---|---|---|---|
| $\lvert D_J \rvert$ | 332180 | 485274 | 635170 | 782753 | 928779 |
| $\lvert P_J \rvert$ | 872120 | 1275229 | 1670157 | 2059566 | 2445064 |

| $J$ | $3,5 \times 10^7$ | $4 \times 10^7$ | $4,5 \times 10^7$ | $5 \times 10^7$ |
|---|---|---|---|---|
| $\lvert D_J \rvert$ | 1073173 | 1216687 | 1358924 | 1500452 |
| $\lvert P_J \rvert$ | 2827226 | 3205958 | 3580839 | 3954180 |

| $\boldsymbol{J}$ | $\mathbf{5,5 \times 10^7}$ | $\mathbf{6 \times 10^7}$ | $\mathbf{6,5 \times 10^7}$ | $\mathbf{7 \times 10^7}$ |
|---|---|---|---|---|
| $\lvert D_J \rvert$ | 1640856 | 1780670 | 1919905 | 2058718 |
| $\lvert P_J \rvert$ | 4325678 | 4695662 | 5064262 | 5431332 |

| $\boldsymbol{J}$ | $\mathbf{7,5 \times 10^7}$ | $\mathbf{8 \times 10^7}$ | $\mathbf{8,5 \times 10^7}$ | $\mathbf{9 \times 10^7}$ |
|---|---|---|---|---|
| $\lvert D_J \rvert$ | 2196834 | 2334373 | 2471678 | 2608393 |
| $\lvert P_J \rvert$ | 5797031 | 6160931 | 6523704 | 6885008 |

| $\boldsymbol{J}$ | $\mathbf{9,5 \times 10^7}$ | $\mathbf{1 \times 10^8}$ |
|---|---|---|
| $\lvert D_J \rvert$ | 2744681 | 2880504 |
| $\lvert P_J \rvert$ | 7245398 | 7605407 |

For the linear regression of $\lvert P_J \rvert$ and $\lvert D_J \rvert$, we used Maple's routine **NonlinearFit** with empirical initial values $a$ and $b$ as below:

**NonlinearFit**($a \times n^b$, X, Y, n, initialvalues = [$a$ = .2, $b$ = .9], output = [leastsquaresfunction, residuals])

We get the following regression:

$$\lvert D_J \rvert \simeq 0.08720724726651 \times J^{0.939843231074621}$$

$$R = 0.999999735986735$$

$$\lvert P_J \rvert \simeq 0.220156649338003 \times J^{0.942278797437848}$$

$$R = 0\,.999999757283397$$